\documentclass[12pt,twoside]{article}
\usepackage{amsmath,amssymb}
\usepackage{graphicx}

\numberwithin{equation}{section}

\title{On the Fractional Mean Value} %
\author{ Abdolali Neamaty,  Mehdi Nategh\\
Department of Mathematics, Mazandaran University, Babolsar, Iran\\
namaty@umz.ac.ir, m.nategh@stu.umz.ac.ir
}
\date{}

\begin{document}
\maketitle
\begin{center}
{\bf\large Abstract}
\end{center}
\par This work, dealt with the classical mean value theorem and took advantage of it in the fractional calculus.
The concept of a fractional critical point
is introduced. Some sufficient
conditions for the existence of a critical point is studied and
an illustrative example relevant to the concept of the time dilation effect is given.
The present paper also includes, some connections between convexity (and monotonicity) with
fractional derivative in the \textit{Riemann}-\textit{Liouville} sense.
\section{Introduction}
\hspace{0.6 cm} Fractional calculus which dates back to more than 300 years ago, has became
one of the thriving areas that is supported by the recent seminal papers together with some special
phenomenological viewpoints. Evidently, there are many dynamic mechanisms in the universe (such
as \textit{anomalous diffusion} and \textit{fractal dynamics}) that are much more better to be described via the fractional calculus than the classical one (see \textit{e.g.} $\cite{JfIp, V E Tarasov, L vlahos}$).
\par
In the present work, the concept of the mean value is nothing more than one is given in the classical calculus. However, the prefix fractional, is to emboss its concept for the next purposes which will be given in sections 2-4.
Also it has been several works in the fractional context that are devoted to this concept.
The generalized mean value theorem was first given in
$\cite{Ijt}$ and this is about a generalized $Taylor$ formula in
a fractional case. In $\cite{Zmo-Ntsh}$, the similar mean value
theorem is obtained but it is more likely to the classical case. This theorem
in both mentioned forms have been used in many papers. In
$\cite{JfIp}$, the \textit{Cauchy} type mean
value is obtained. Moreover, the mean value is studied in a more
abstract and operational form. Another generalized mean value
theorem is proved in $\cite{Kdie}$ with slight changes, compared to the one is
obtained in $\cite{Ijt}$.
\par The accurate relation between convexity and increase of the differential in local sense, depends on the usual order on the real line. Connecting the convexity to the differential
(in non-local sense) in the same manner, needs bringing an order on the domains of integrations, since we are
dealing with the integro-differential operators. This is done in section 4. Moreover, a regularization is introduced, that upon which,
the distribution of the fractional mean values have became uniform. This is because we make use of the fractional mean value
in \textit{Theorem 4.5.}
\par This paper is organized as follows:
In section 2, some preliminaries on the concept of the mean value
is introduced and an approximation of it as a root of a
polynomial is given.\\
In section 2, the concept of the fractional critical point as the root
of the fractional derivative is introduced, and sufficient conditions
for its existence are studied. Also, an interpretation of the
fractional critical point is given based on the time dilation
effect (according to $\cite{Ip2}$).\\
Connections between fractional derivative and convexity (as mentioned above) and
monotonicity are studied in section 3. A fundamental condition for
this connection is introduced. This condition, defines some order on
the fractional derivatives with respect to the orders on the intervals
that are being involved in the \textit{Riemann-Liouville} integro-differential operator.\\
\section{Fractional mean value}
\textbf{Definition 2.1.} $\cite{Ak}$ Let $[a,b]$ be an
interval and $f$ be an integrable function. The left
\textit{Riemann}-\textit{Liouville\textit} fractional integral and
derivative of order
$\alpha \in (0,1)$, (provided they exist) are defined by\\
\begin{eqnarray}
I^{\alpha}_{a^{+}}f(x): &=& \frac {1}{\Gamma
(\alpha)}\int_{a}^{x}\frac
{f(t)}{(x-t)^{1-\alpha}}dt\,\,\,\,\,\,\,\forall x \in (a,b],\\
D^{\alpha}_{a^{+}}f(x) : &=& \frac {d}{dx} I^{ 1- \alpha
}_{a^{+}}f(x) =\frac {1}{\Gamma ( 1 - \alpha)} \frac {d}{dx}
\int_{a}^{x}\frac
{f(t)}{(x-t)^{\alpha}}dt\,\,\,\,\,\,\,\forall x \in (a,b].
\end{eqnarray}
where $\Gamma(\alpha)$ is the Gamma function. The space $I^{\alpha}_{a^{+}}(L^{p})$ for $1 \leq p \leq \infty$, is
defined by
\begin{equation}
I^{\alpha}_{a^{+}}(L^{p}):=\Big\{f \,\, \Big |\,\,
f=I^{\alpha}_{a^{+}}\phi,\,\, \phi \in L^{p}(a,b)\Big\}.
\end{equation}
The left \textit{Caputo} fractional derivative of a function $f$
(provided it exists) is defined by
\begin{equation}
^{C}D^{\alpha}_{a^{+}}f(x):=\frac {1}{\Gamma
(1-\alpha)}\int_{a}^{x}\frac
{f'(t)}{(x-t)^{\alpha}}dt\,\,\,\,\,\,\,\forall x \in (a,b].
\end{equation}
The set of absolutely continuous functions, that is denoted by
$AC[a,b]$, is the set of all functions $f$, which has representation
of the form ($\cite{Ak}$):
\begin{equation}
f(x)=c + I^{n}_{a^{+}}\phi, \,\,\,\,\,\,\,\,  c \in \mathbb{R}, \,\,
\phi \in L( a,b ).
\end{equation}
\textbf{ Proposition 2.2. } $\cite{Ak}$
\begin{description}
    \item  {a)} If $f \in L^{p}(a,b)$, ($ 1 \leq p < +\infty$), then
       \begin{equation}
      D^{\alpha}_{a^{+}} I^{\alpha}_{a^{+}} f(x) = f(x).
      \end{equation}
    \item  {b)} If $f\in L^{1}(a,b)$  and $f_{1-\alpha}=I^{1-\alpha}_{a^{+}}f \in
AC[a,b]$, then
\begin{equation}
 I^{\alpha}_{a^{+}}D^{\alpha}_{a^{+}}f(x) = f(x)- \frac
{f_{1-\alpha}(a)}{\Gamma(\alpha)}(x-a)^{\alpha -1},
\end{equation}
holds almost everywhere on $[a,b]$.
     \item  {c)} If $f\in I^{\alpha}_{a^{+}}(L^{p})$, $(1\leq p \leq \infty)$
            then
\begin{equation}
I^{\alpha}_{a^{+}}D^{\alpha}_{a^{+}}f(x)=f(x).
\end{equation}
\end{description}
\textbf{ Proposition 2.3.} Assume that $ f $ is differentiable on $[ a,b ]$ with $ f(a)=0 $, then\\
$(a)$ $\cite{Ak}$,
\begin{equation}
D^{\alpha}_{a^{+}}f(x)=\,
 ^{C}D^{\alpha}_{a^{+}}f(x)=I^{1-\alpha}_{a^{+}}f'(x).
\end{equation}\\
$(b)$  $\cite{Ip}$,
\begin{equation}
\lim_{ \alpha \rightarrow 0^{+} } I^{\alpha}_{a^{+}} f(x) = f(x).
\end{equation}
\textbf{ Proposition 2.4.} $\cite{S-s A-k}$\\
$(a)$ If $f \in AC[ a,b ]$, then $f_{1-\alpha} \in AC [ a,b ]$.\\
$(b)$ \textit{Abel}'s integral equation is solvable in $L^{1}(a,b)$
if and only if $f_{1-\alpha} \in AC [a,b]$ and we have
\begin{equation}
f_{ 1-\alpha }(x) = \frac {1}{ \Gamma(2-\alpha) } \Big \{ f(a)
(x-a)^{1-\alpha} + \int_{a}^{x}f'(t)()^{x-t1-\alpha} dt \Big \}.
\end{equation}\\\\
\textbf{ Assumptions 2.5. } We assume that $f \in C^{1}[a,b]$
with $f(a)=0$ and $D^{\alpha}_{a^{}+} f \in C [ a,b ]$.\\\\
In this paper, all statements except
\textit{Proposition 4.6}, presume \textit{Assumption 2.5}. Also
\textit{ Proposition 2.4}, together with \textit{Assumption 2.5}, yield the
existence of $D^{\alpha}_{a^{+}}f(x)$ for every $\alpha \in ( 0,1
)$, since $C^{1}[a,b] \subset AC [a,b]$.

Now we begin with this known lemma in the first semester calculus:
Let $\phi :\in L^{1}(a,b)$, then
\begin{equation}
\exists \,\xi \in
(a,b),\,\,\,\,f(\xi).\int_{a}^{b}\phi(t)dt=\int_{a}^{b}\phi(t)f(t)dt.
\end{equation}
Letting $\phi_{x}(t)=\frac {(x-t)^{-\alpha}}{\Gamma(1-\alpha)}$ we
have
\begin{equation}
\exists \,\xi \in (a,x),\,\,\,\,f(\xi).\frac
{(x-a)^{1-\alpha}}{\Gamma(2-\alpha)} =\frac
{1}{\Gamma(1-\alpha)}\int_{a}^{x}\frac
{f(t)}{(x-t)^{\alpha}}dt=I^{1-\alpha}_{a^{+}}f(x).
\end{equation}
We name $\xi$ the \textbf{fractional mean value}. Since $\xi$ is not
unique, we define
\begin{equation}
\xi(x,\alpha):= \sup \Big\{ \xi \in (a,x)\,|\,
f(\xi)=\Gamma(2-\alpha).\Big(I^{1-\alpha}_{a^{+}}f(x)\Big).(x-a)^{\alpha-1}
  \Big\}.
\end{equation}
In general, all possible $\xi$'s are belonging to the set
$\Lambda(\alpha,f,x)$ which is defined by
\begin{equation}
\Lambda(\alpha,f,x):=f^{-1}\Big\{
\Gamma(2-\alpha).\Big(I^{1-\alpha}_{a^{+}}f(x)\Big).(x-a)^{\alpha-1}
\Big\}.
\end{equation}
The set $\Lambda$ is closed, since $f$ is continuous. Let $ g(x):=
\Gamma(2-\alpha).\Big(I^{1-\alpha}_{a^{+}}f(x)\Big).(x-a)^{\alpha-1}$.
If $f$ is monotone on $( a,b )$, then for $x \in (a,b)$, there
exists a unique $\xi ( x,\alpha ) \in ( a,x )$ for which $ f(
\xi(x,\alpha ) ) = g(x) $. In other words, $\forall x \in (a,b)$,
$g(x) \in Im ( f ).$
\\\\
\textbf{ Proposition 2.6. } Suppose $f$ is monotone on $(
a,b )$. Then
$\xi(x,\alpha) \in C^{1}( a,b )$.\\

The following lemma gives some details about the mean value function
$\xi$, with respect to the variable $x$, but under the certain
condition.\\\\
\textbf{ Proposition 2.7.} Suppose that $h:[a,b] \rightarrow
[a,b]$ is differentiable. Assume that there exists some $x_{0} \in
(a,b)$ at which $h$ has a local extremum and also
$h(x_{0})=\xi(x_{0})$, where $\xi(x) \in \Lambda(\alpha,f,x)$. Then
there exists some $x \in (a,b)$, such that
\begin{equation}
D^{\alpha}_{a^{+}}f(x)=D\Big(f(h(x)).\frac
{(x-a)^{1-\alpha}}{\Gamma(2-\alpha)}\Big).
\end{equation}\\
\textbf{Proof.} In contrary assume that for every $x \in
(a,b)$, we have:
\begin{equation}
D^{\alpha}_{a^{+}}f(x) \neq D\Big(f(h(x)).\frac
{(x-a)^{1-\alpha}}{\Gamma(2-\alpha)}\Big).
\end{equation}
Letting
$$
F(x,h)=I^{1-\alpha}_{a^{+}}f(x)-f(h(x)).\frac
{(x-a)^{1-\alpha}}{\Gamma(2-\alpha)},
$$
we have $F(x_{0},h(x_{0})=0$ and $\frac {\partial F}{\partial
x}(x_{0},h(x_{0})) \neq 0$. By Implicit Function theorem, $x$ will
be a function of $h$ in the vicinity of $x_{0}$ and this is
impossible since $h'(x_{0}) = 0$. This proves the assertion of the
Lemma.
$\Box$\\

Regarding to \textit{Proposition 2.7}, when a given function $h \in
C^{1}([a,b],[a,b])$, satisfies $(2.17)$ everywhere, one conclude
that $h$
is monotone, or $h(x) \neq \xi(x)$ holds everywhere.\\\\
\textbf{ Proposition 2.8.} Let $f \in C^{n+1}[a,b]$. Then
$\xi(\alpha)$ can be estimated by finding the roots of the following
polynomial
\begin{eqnarray}
\nonumber &\forall & \delta <b-a, \\ \nonumber &&
\sum_{j=1}^{n}f^{(j)}(a).\Big\{ \frac
{\delta^{1-\alpha}}{\Gamma(2-\alpha)j!}x^{j}-\frac
{\delta^{j+1-\alpha}}{\Gamma(j+2-\alpha)} \Big\} \\ &\,\,&-
\Big[I^{n+2-\alpha}_{a^{+}}f^{(n+1)}\Big](a + \delta) + o(\frac
{1}{n!})=0,
\end{eqnarray}
where $x=\xi -a$.\\\\
\textbf{Proof.} First observe that:
\begin{equation}
f(\xi).\frac {\delta^{1-\alpha}}{\Gamma(2-\alpha)}
=\Big[I^{1-\alpha}_{a^{+}} f \Big](a + \delta)=
\sum_{j=0}^{n}f^{(j)}.\frac
{\delta^{j+1-\alpha}}{\Gamma(j+2-\alpha)} +
\Big[I^{n+2-\alpha}_{a^{+}}f^{(n+1)}\Big](a + \delta).
\end{equation}
On the other hand, taking into account the expansion of $f(\xi)$
about the point $a$ that is:
\begin{equation}
f(\xi)=\sum_{j=0}^{n}f^{(j)}(a).\frac {(\xi-a)^{j}}{j!} + o(\frac
{1}{n!}),
\end{equation}
and equalizing with $f(\xi)$ in $(2.19)$, gives the result. $\Box$\\
\section{Fractional critical point}
\hspace{0.6 cm} In this section, the roots of the fractional
derivative of a function $f$, will be studied. This may be thought
of as a parallel study to the critical points in the ordinary calculus.
For the similarity, we name the root of $D^{\alpha}_{a^{+}}f$, the
\textit{fractional critical point} of $f$.
\par To find the fractional critical point of a function $f$ of order
$1-\alpha$, it is necessary and sufficient to find the nodes, those
are intersections of two curves $(x-a)^{\alpha}$ and $f(\xi(x))$.
This can be stated in a rigorous way as follows:\\\\
\textbf{ Lemma 3.1.} Suppose $f$ is monotone, then
\begin{equation}
f(\xi (x))=(x - a)^{\alpha}, \,\,\,\,\,\,\,\,\,  \xi (x) \in \Lambda ( \alpha , f , x )
\end{equation}
if and only if $D^{1-\alpha}_{a^{+}}f(x)=0$.\\\\
\textbf{Proof.} Let $g(x)=f( -\eta(x))$, where $\eta(x)= - \xi (x) $. By \textit{Proposition 2.6}, $\eta$ is
differentiable. The point $x_{0} \in ( a,b )$ satisfies the equation
$g(x_{0})=(x_{0}-a)^{\alpha}$, if and only if it satisfies the
following equation
\begin{equation}
(x-a)g'(x) - \alpha g(x)=0.
\end{equation}
Multiplying $(3.2)$ by $\frac {(x-a)^{\alpha -1}}{\Gamma(1+\alpha)}$
we get:
\begin{eqnarray}
\nonumber 0 &=& \frac {1}{\Gamma(1+\alpha)} \Big( -
(x_{0}-a)^{\alpha}f'( - \eta(x_{0}))\eta'(x_{0}) -
\alpha(x_{0}-a)^{\alpha -1}f( -\eta(x_{0}))\Big)\\ \nonumber &=& -
\frac {d}{dx} \Big(\frac
{(x-a)^{\alpha}}{\Gamma(1+\alpha)}.f( \xi(x))\Big)_{x=x_{0}} \\
\nonumber &=& - \frac {d}{dx}I^{\alpha}_{a^{+}}f(x)_{\Big|_{x=x_{0}}} \\
&=& - D^{1-\alpha}_{a^{+}}f(x_{0}).\,\,\,\,\, \Box
\end{eqnarray}\\
\textbf{ Theorem 3.2.} Assume that $f(x)=0$ for some $x \in (a,b)$, then there
exists some $\xi \in
(a,x]$ for which $D^{\alpha}_{a^{+}}f(\xi)=0$.\\\\
\textbf{  Proof.} By $(2.11)$,  $f_{1-\alpha}(a)=0$. Also
by \textit{Proposition 2.2 (b)}, we have
\begin{equation}
\exists \, \xi \in \Lambda(\alpha,D^{\alpha}_{a^{+}}f,x), \,\,
D^{\alpha}_{a^{+}}f(\xi(x)).\frac
{(x-a)^{\alpha}}{\Gamma(\alpha)}=I^{\alpha}_{a^{+}}D^{\alpha}_{a^{+}}f(x)=f(x)=0,
\end{equation} $\Box$\\
\textbf{ Remark 3.3:} According to the interpretation of
the fractional derivative, that is introduced in $\cite{Ip2}$, let
$v(t)$ be the velocity measured by an individual, that is trapped in
a gravitational field, $V$ be the velocity measured by an
independent observer, and $\tau$ and $t$ be the time in the presence of
the gravitational field and the cosmic time respectively, then we
have:
\begin{equation}
V(t)=D^{\alpha}_{a^{+}}v(t)=\frac {1}{\Gamma(1-\alpha)}\frac
{d}{dt}\int_{a}^{t}\frac {v(\tau)}{(t-\tau)^{\alpha}}d \tau.
\end{equation}
If $v(t)$ vanishes at specific time $t$, then \textit{Theorem 3.2}
states that, the independent observer has measured zero velocity, sooner
at some time $\xi(t) \leq t$.\\\\
\textbf{ Lemma 3.4.} For $x \in [ a,b ]$, we have
\begin{equation}
\left\{\begin{array}{ll} \lim_{\alpha \rightarrow
0^{+}}D^{\alpha}_{a^{+}}f(x)=f(x),
\\ \lim_{\alpha \rightarrow
1^{-}}D^{\alpha}_{a^{+}}f(x)=f'(x).
\end{array}\right.
\end{equation}
\textbf{ Proof.} By assumption we have (\textit{Proposition}
$2.3 (a)$)
\begin{equation}
D^{\alpha}_{a^{+}}f(x)\,=\,^{C}D^{\alpha}_{a^{+}}f(x)=
I^{p}_{a^{+}}f'(x)
\end{equation}
where $p=1-\alpha$. Now by \textit{Proposition} $2.3 (b)$
\begin{equation}
\lim_{\alpha \rightarrow 1^{-}}D^{\alpha}_{a^{+}}f(x)=
\lim_{p\rightarrow 0^{+}}I^{p}_{a^{+}}f'(x)= =f'(x).
\end{equation}
Again by $D^{\alpha}_{a^{+}}f(x)=I^{1-\alpha}_{a^{+}}f'(x)$ we have
\begin{eqnarray}
\nonumber \Big|I^{1-\alpha}_{a^{+}}f'(x)-I^{1}_{a^{+}}f'(x)\Big|
&\leq& M\Big|\int_{a}^{x}\Big(\frac
{(x-t)^{-\alpha}}{\Gamma(1-\alpha)}-1\Big)dt\Big|\\
&=& M\Big| \frac {(x-a)^{1-\alpha}}{\Gamma(2-\alpha)}-(x-a)\Big|.
\end{eqnarray}
where $M= \max_{t \in [a,b]} |f'(t)|$. Taking the limit $\alpha
\rightarrow 0^{+}$, gives the result.
$\Box$\\

Let $ X_{\epsilon, \alpha}:=\Big \{ r \in B(x_{0}, \epsilon)  \, | \,
D^{\alpha}_{a^{+}} f( r ) = 0 \Big \} $ and $r(\alpha):= \sup X_{\epsilon, \alpha}$ ( provided
$ X_{\epsilon, \alpha} \neq \emptyset $ ), where $ \epsilon \in \mathbb{R}^{+} $, $ \alpha \in [ 0,1 ) $ and $ x_{0} \in [ a,b ] $. We have the following statement:\\\\
\textbf{ Theorem 3.5.} Suppose $x_{0}$ and $x_{1} $ are the
unique maximum (minimum) and root of $f$ in $(a,b]$
respectively, then
\begin{equation}
\left\{\begin{array}{ll} \lim_{\alpha \rightarrow 0^{+}}
r(\alpha)=x_{1},
\\ \lim_{\alpha \rightarrow 1^{-}} r(\alpha)= x_{0}.
\end{array}\right.
\end{equation}
\textbf{ Proof.} Without loss of generality, assume the
$x_{0}$ is a maximum. For $t \in (a,b) \backslash \{x_{0}\}$ we have
\begin{equation}
\lim_{\alpha \rightarrow 1^{-}}D^{\alpha}_{a^{+}}f(t)=
\left\{\begin{array}{ll} f'(t) >0 & t \in (a,x_{0}),
\\ f'(t)<0 & t \in (x_{0},b).
\end{array}\right.
\end{equation}
Since for every $\epsilon \in \mathbb{R}^{+}$, $f'$ changes sign
over $B(x_{0},\epsilon)$, for $\alpha$ close enough to 1,
$D^{\alpha}_{a^{+}}f$ changes sign too and therefore, $X_{\epsilon, \alpha}$
is nonempty. Let $ r =\limsup_{\alpha \rightarrow 1^{-}}
r(\alpha)$. If $r \neq x_{0}$, then by continuity of $f'$ and
$ \lim_{ \alpha \rightarrow 1^{-} }D^{\alpha}_{a^{+}} f (r( \alpha )) = f' ( r ( \alpha ) )$, one obtain
\begin{eqnarray}
0 = \limsup_{ \alpha \rightarrow 1^{-} }
D^{\alpha}_{a^{+}}f( r(\alpha))  = f' ( r ) \neq 0.
\end{eqnarray}
Similarly we infer that, $\liminf_{\alpha
\rightarrow 1^{-}} r(\alpha)=x_{0}$. The same result can be obtained in the case
$ \alpha \rightarrow 0^{+} $. $\Box$\\

\section{Convexity}
\hspace{0.6 cm} In classical derivative, convexity of a
differentiable function is equivalent to that function has
a nondecreasing derivative. In this section, our aim is to give a
similar argument. But first we need an order to give a meaning to
the concept of increase or decrease of the fractional derivative.
\par Regarding to $\cite{Ip2}$, upon which, the fractional
derivative, interpreted as the real velocity (with respect to cosmic
time) is calculated for an event that happens in the presence of some
strong gravitational field (with time dilation effect), we propose
the following:

With respect to the non-locality of the fractional derivative that
involves integrating over a time interval namely $[a,x]$, it seems
to be reasonable if we define an order on the fractional derivative,
in such a manner that we fix a constant length for two different
time domains that the integrations are taken over them, that is:\\\\
\textbf{ Definition 4.1.} A real valued function $f$ is called \textit{$\delta$-increasing}, if
for every $x_{0},y_{0}, x,y \in
\mathbb{R}$ and $ \delta > 0$ with $ x_{0} < x < y_{0} < y $ and $x - x_{0} = y - y_{0} = \delta$
we have
$$
D^{\alpha}_{x_{0}^{+}}f(x - x_{0}) \leq D^{\alpha}_{y_{0}^{+}}f(y -y_{0}).
$$
Next step is to set a discipline for the mean value corresponding to
the fractional integration of order $1-\alpha$:\\\\
\textbf{ Definition 4.2.} We say that a continuous function
$f$ has the property $(\textbf{P})$, if for every $x_{0},y_{0}, x, y \in
\mathbb{R}$ and $\delta > 0$ with $ x_{0} < x < y_{0} < y $ and $ x - x_{0} = y - y_{0} = \delta $, we have
$$
\xi_{x} - \xi_{y}= y_{0} - x_{0},
$$
where $\xi_{x_{0}+ \delta} \in \Lambda(\alpha,f \circ ( x - x_{0} )  ,x_{0}+ \delta)$ and
$\xi_{y_{0}+ \delta} \in \Lambda(\alpha,f \circ ( y - y_{0} ) ,y_{0}+ \delta)$.\\

The property $(\textbf{P})$ is about the independence of the mean
value $\xi$ to the initial point $x_{0}$ or $y_{0}$, when we are
dealing with the fractional differentiation that involves
integrating and the integration is taken over two different time
intervals of the same length,
$[x_{0},x_{0} + \delta]$ and $[y_{0},y_{0} + \delta]$.\\\\
\textbf{ Example 4.3. } Let $ f(x)=x^{\beta} $ ( $ \beta \in \mathbb{R} $ ), then $f$ is \textit{$\delta$-increasing}
since
\begin{equation}
D^{\alpha}_{x_{0}^{+}} f(x - x_{0}) = \frac { \Gamma (1 + \beta) }{ \Gamma( \beta - \alpha +1 ) }
( x - x_{0} )^{ \beta - \alpha } = \frac { \Gamma (1 + \beta) }{ \Gamma( \beta - \alpha + 1 ) }
\delta^{ \beta - \alpha }.
\end{equation}
Also we have
\begin{equation}
I^{1 - \alpha}_{x_{0}^{+}} f(x - x_{0})^{\beta} = \frac { \Gamma ( \beta ) }{ \Gamma (\beta - \alpha + 1) }
( x - x_{0} )^{ \beta - \alpha + 1 } = \frac { ( x - x_{0} )^{ 1 - \alpha } }{ \Gamma ( 2 - \alpha ) } ( \xi (x) - -x_{0} )^{ \beta },
\end{equation}
and we obtain
\begin{equation}
\frac { \xi (x) - x_{0} }{ x - x_{0} } = \Bigg \{ \frac { \Gamma ( 2 - \alpha ) \Gamma ( \beta ) }{ \Gamma ( \beta - \alpha + 1 ) } \Bigg \}^{ \frac {1}{ \beta }}.
\end{equation}
For $x_{0} < x < y_{0} < y $ with $ x - x_{0} = y - y_{0} = \delta $, we have $ \xi (x) - x_{0}  = \xi (y) - y_{0} $ if and only if $ \Bigg \{ \frac { \Gamma ( 2 - \alpha ) \Gamma ( \beta ) }{ \Gamma ( \beta - \alpha + 1 ) } \Bigg \}^{ \frac {1}{ \beta }} $ is a positive real number and it is equivalent to $f$ satisfies the property \textbf{$ (P) $}.\\\\
\textbf{ Remark 4.4.} Under the monotonicity assumption of
$f$ with property $(\textbf{P})$ and denoting $\xi(x,a)$ the mean
value defined by $(2.14)$, we have
\begin{equation}
\frac {\partial}{\partial a}\Bigg\{
\Gamma(2-\alpha)\Big(I^{1-\alpha}_{a^{+}}f\Big)(x).(x-a)^{\alpha
-1}\Bigg\} = \frac {\partial}{\partial a}f(\xi (x,a)).
\end{equation}
Indeed because of the independence of $\xi(x)-a$ to the initial
point $a$, there must be a function $h(x)$ such that
$\xi(x)=a+h(x)$. So we have
$$
h(x)=f^{-1}\Bigg\{
\Gamma(2-\alpha)\Big(I^{1-\alpha}_{a^{+}}f\Big)(x).(x-a)^{\alpha-1}
\Bigg\}-a.
$$
Now differentiating with respect to $a$ and utilizing
$(f^{-1})'=\frac {1}{f'}$, gives $(4.1)$.\\\\
\textbf{ Theorem 4.5.} Let $f'$ satisfies the
property $(\textbf{P})$. Then $f$ is convex if and only if it is \textit{$\delta$-increasing}.\\\\
\textbf{ Proof.} For $x_{0},y_{0} \in [a,b]$ and positive $\delta$ with
$x_{0} < x < y_{0} < y $ and $ x- x_{0} = y - y_{0} = \delta$, we have
\begin{eqnarray}
\nonumber  D^{\alpha}_{x_{0}^{+}}f(x - x_{0}) - D^{\alpha}_{y_{0}^{+}}f(y - y_{0}) & = &
I^{1-\alpha}_{x_{0}^{+}}f'(x - x_{0}) - I^{1-\alpha}_{y_{0}^{+}}f'(y - y_{0}) \\ &=&
\frac {\delta^{1-\alpha}}{\Gamma(2-\alpha)}. \Bigg(f'(\xi(x)) -
f'(\xi(y))\Bigg).
\end{eqnarray}
The fractional derivative is increasing if and only if $f'(\xi(x))
\leq f'(\xi(y))$ and by continuity of $f'$, we infer that $f'$ is
increasing and therefore $f$ is convex. $\Box$\\

For a positive $\tau$, let $\Delta f(x):=f(x + \tau)-f(x)$. If for a
constant $\tau$, we have $\Delta f(x) \geq 0$, then we say that $f$ is
increasing up to the positive constant $\tau$. Next Lemma is about the
relationship between the monotonicity of a function and its
fractional derivative.\\\\
\textbf{ Proposition 4.6.} Let $f \in C[0,b]$ and $f(x)=0$ for $x
\in \mathbb{R} \backslash [0,b]$. Let $D^{\alpha}_{a^{+}}f$ exists
and for $0 < \tau < b$, we have $D^{\alpha}_{a^{+}}f(x) \leq
D^{\alpha}_{a^{+}}f(x + \tau)$.
Then $f$ is increasing (up to $\tau$) if $\Delta f(0) \geq 0$.\\\\
\textbf{ Proof.} The difference $\Delta
D^{\alpha}_{0^{+}}f(x)$ (that is a function of $x$ and $\tau$) is
positive. Taking the \textit{Laplace} transform we obtain
\begin{equation}
L \Delta D^{\alpha}_{0^{+}}f(x)=s^{\alpha} \Big(Lf(x + \tau) - Lf(x)\Big)
- \Delta f(0).
\end{equation}
Then taking the inverse gives
\begin{eqnarray}
\nonumber \Delta f(x) &=&  \frac {1}{ \Gamma (\alpha) } \Big (
x^{\alpha -1}\Delta f(0) + x^{\alpha -1}
* \Delta D^{\alpha}_{a^{+}}f(x) \Big ) \\ & = &\frac { x^{\alpha -1} }{ \Gamma (\alpha) } \Delta
f(0) + I^{\alpha}_{a^{+}}\Delta D^{\alpha}_{a^{+}}f(x) \geq
0.,
\end{eqnarray}
where asterisk means convolution. $\Box$\\\\
\textbf{ Remark 4.7.}\\
It is easy to verify that, if $D^{\alpha}_{a^{+}}f(x) \leq
D^{\alpha}_{a^{+}}g(x)$ and $f(0)=g(0)$, then $f(x) \leq g(x)$.\\
Let $f$ be periodic with period $\tau$. Then $(4.7)$ gives
$I^{\alpha}_{a^{+}}\Delta D^{\alpha}_{a^{+}}f(x)=0$ and
therefore
\begin{equation}
D^{\alpha}_{a^{+}}f(x + \tau)=D^{\alpha}_{a^{+}}f(x).
\end{equation}
In other words, periodic functions have periodic fractional
derivatives with the same period.\\\\
\textbf{Conclusion} In this paper, the fractional mean value is introduced and its approximation as the roots of a polynomial is obtained.
Fractional critical point is introduced, some sufficient conditions for its existence are studied and an interpretation of such critical point
is given with respect to the time dilation effect in a gravitational field. Moreover, Convexity and monotonicity and their possible connections to
increase of a fractional derivative in the \textit{Riemann-Liouville} sense are discussed.

\end{document}